\documentclass {article}
\usepackage{authblk}
\usepackage{etex}
\usepackage[utf8]{inputenc}
\usepackage[english]{babel}
\usepackage{amsmath}
\usepackage{amsthm}
\usepackage{amssymb}
\usepackage{sectsty}
\usepackage{titlesec}
\usepackage{color}
\usepackage[color,matrix,arrow]{xy}
\usepackage{amsgen}
\usepackage{amstext}
\usepackage{amsbsy}
\usepackage{amsopn}
\usepackage{amsfonts}
\usepackage{eepic}
\usepackage{graphicx}
\usepackage{epsf}
\usepackage{pstricks}
\usepackage{enumerate}
\usepackage{eqnarray}
\usepackage{faktor}
\xyoption{all}

\usepackage{pgf}
\usepackage{tikz-cd}
\usetikzlibrary{automata, arrows.meta, positioning,decorations.pathmorphing}

\newcommand{\footrecall}[1]{%
} 
\usetikzlibrary{snakes,shapes,arrows,automata}

\titleformat*{\section}{\large\bfseries}
\titleformat*{\subsection}{\normalsize \bfseries}

\newcommand{\N}{\mathbb{N}}
\newcommand{\Z}{\mathbb{Z}}

\newcommand{\Q}{\mathbb{Q}}
\newcommand{\Fix}{\text{Fix}}

\newcommand{\End}{\text{End}}
\newcommand{\Ker}{\text{Ker}}

\newcommand{\Inn}{\text{Inn}}

\newcommand{\Orb}{\text{Orb}}
\newcommand{\Aut}{\text{Aut}}
\newcommand{\Via}{\text{Via}}

\newcommand{\Rat}{\text{Rat}}
\DeclareMathOperator*{\bigplus}{\scalerel*{+}{\textstyle\sum}}
\usepackage{scalerel}
\newcommand{\mc}{\mathcal}

\theoremstyle{definition}
\newtheorem{theorem}{Theorem}[section]
\newtheorem{corollary}[theorem]{Corollary}

\newtheorem{proposition}[theorem]{Proposition}

\newtheorem{lemma}[theorem]{Lemma}
\newtheorem{example}[theorem]{Example}

\begin{document}
 
 
\title{On generalized conjugacy and some related problems}
\author{Andr\'e Carvalho} 
\maketitle

\begin{abstract}
We establish a connection between the generalized conjugacy problem for a $G$-by-$\Z$ group, $GCP(G\rtimes \Z)$, and two algorithmic problems for $G$: the generalized Brinkmann's conjugacy problem, $GBrCP(G)$, and the generalized twisted conjugacy problem, $GTCP(G)$. We explore this connection for generalizations of different kinds: relative to finitely generated subgroups, to theirs cosets, or to recognizable, rational, context-free or algebraic subsets of the group. Using this result, we are able to prove that $GBrCP(G)$ is decidable when $G$ is a virtually polycyclic group, which implies in particular that the generalized Brinkmann's equality problem, $GBrP(G)$, is decidable if $G$ is a finitely generated abelian group. Finally, we prove that if $G$ is a finitely generated virtually free group, then the simple versions of Brinkmann's equality problem and of the twisted conjugacy problem, $BrP(G)$ and $TCP(G)$, are decidable.
\end{abstract}

\section{Introduction}
The conjugacy problem on $G$, $CP(G)$,  consists on deciding whether two elements, given as input, are conjugate or not in $G$. 

The twisted conjugacy problem on $G$, which we denote by $TCP(G)$ consists of deciding  whether, taking two elements $g,h\in G$ and an automorphism $\phi\in Aut(G)$, there is some $x\in G$ such that $h=(x^{-1}\phi)gx$.
In this case, we say that $g$ and $h$ are $\phi$-twisted conjugate. 

In \cite{[Bri10]}, Brinkmann proved that, in case $G$ is a free group of finite rank, there is an algorithm to decide whether, on input $g,h\in G$ and $\phi\in Aut(G)$, there is some $k\in \Z$ such that $g\phi^k$ is conjugate to $h$. We call this problem Brinkmann's Conjugacy Problem and denote it by $BrCP(G)$. Brinkmann also solved a similar problem with equality instead of conjugation in  \cite{[Bri10]}. We will call this Brinkmann's (Equality) Problem and denote it by $BrP(G)$. 
This notion was later generalized in \cite{[BMV10]} by the notion of orbit decidability. For more on this, the reader is referred to the survey \cite{[Ven14]}.

In this paper, we propose the following generalizations of these problems: 
let $\mc C$ be a class of subsets of $G$.
\begin{itemize}
\item \textbf{$GCP_{\mc C}(G)$ -  $\mc C$-generalized conjugacy problem:} taking as input a subset $K\in \mc C$ and an element $x\in G$, decide whether there is a conjugate of $x$ in $K$;
\item \textbf{$GTCP_{\mc C}(G)$ - $\mc C$-generalized twisted conjugacy problem:} taking as input a subset $K\in \mc C$, an automorphism $\phi\in \Aut(G)$ and  an element $x\in G$, decide whether $x$ has a $\phi$-twisted conjugate  in $K$, i.e., whether there is some $z\in G$ such that $z^{-1}\phi x z\in K$;
\item \textbf{$GBrCP_{\mc C}(G)$  -  $\mc C$-generalized Brinkmann's conjugacy problem:} taking as input a subset $K\in \mc C$, an automorphism $\phi\in \Aut(G)$ and  an element $x\in G$, decide whether there is some $k\in \Z$ such that $x\phi^k$ has a conjugate in $K$;
\item \textbf{$GBrP_{\mc C}(G)$ -  $\mc C$-generalized Brinkmann's problem:} taking as input a subset $K\in \mc C$, an automorphism $\phi\in \Aut(G)$ and an element $x\in G$, decide whether there is some $k\in \Z$ such that $x\phi^k$ belongs to $K$.
\end{itemize}

Some natural classes to consider are $f.g.$,  $Rat$, $Rec$, $Alg$ and $CF$, the classes of finitely generated subgroups, rational, recognizable, algebraic and context-free subsets, respectively (for more information on these classes the reader is referred to \cite{[Ber79],[BS21],[Car22b]}). A seemingly less natural class to consider is the class of cosets of finitely generated groups, which we will denote by $[f.g \, coset]$. As will become clear, this  will be more adequate to us than the class of $f.g.$ 

We have the following inclusions:
\[
\begin{tikzcd}  
f.g \arrow[r,  phantom, sloped, "\subseteq"] & f.g\, coset  \arrow[r,  phantom, sloped, "\subseteq"] &Rat \arrow[r, phantom, sloped,  "\subseteq"] &Alg\\
&&Rec\arrow[u,  phantom, sloped, "\subseteq"]\arrow[r,  phantom, sloped, "\subseteq"] &CF\arrow[u,  phantom, sloped, "\subseteq"]
\end{tikzcd}.
\]

Whenever we want to consider an instance of an algorithmic problem with restrictions on the input, we will write the restrictions as indices. For example, if $H\leq  G$, $\phi\in \Aut(G)$ and $x,y\in G$, we write $GTCP_{(yH,\phi, x)}(G)$ to denote the problem of deciding whether $x$ has a $\phi$-twisted conjugate in $yH$ or not. Notice that it might happen that not all components of the input (subset, automorphism and element) are restricted, but some are. For instance, we might write $GTCP_{([f.g \, coset],\phi)}(G)$ meaning that we take as input a subgroup $H\leq G$ and elements $x,y\in G$ and we want to decide whether $x$ has a $\phi$-twisted conjugate element in $yH$.  Also, we will use the symbol $\equiv$ to denote equivalence between  two decision problems and $+$ to denote the OR operator, in the sense that if we have two decision problems ${P_1}_{(Input_1)}$ and ${P_2}_{(Input_2)}$, then the answer of ${P_1}_{(Input_1)}+{P_2}_{(Input_2)}$ is \texttt{YES}  if the answer of  ${P_1}_{(Input_1)}$ or the answer of  ${P_2}_{(Input_2)}$ is \texttt{YES} and  \texttt{NO} if the answer of both   ${P_1}_{(Input_1)}$ and  ${P_2}_{(Input_2)}$ is \texttt{NO}.

We remark that not much is known about these generalized problems. In \cite{[LS11]},  Ladra and Silva work in the context of rational subsets and solve $GCP_{Rat}(G)$ when $G$ is a  finitely generated virtually free group. They do so by solving $GTCP_{Rat, \phi}(F_n)$, for the very particular case where  $\phi$  is a virtually inner automorphism, i.e., there is a power of $\phi$ that belongs to $\Inn(F_n)$. In the survey \cite{[Ven14]}, Ventura remarks that not much is known about $GBrP_{f.g.}(G)$ even when $G$ is a free  or a free-abelian group.

Let $A=\{a_1,\ldots, a_n\}$, $G=\langle A\mid R\rangle$ be a group and $\phi\in \Aut(G)$. A $G$-by-$\Z$ group has the form 
\begin{align}
\label{presentation semidirect}
G\rtimes_\phi \Z=\langle A,t \mid R, t^{-1}a_it= a_i\phi\rangle.
\end{align}

Every element of $G\rtimes_\phi \Z$ can be written in a unique way as an element of the form $t^ag$, where $a\in \Z$ and $g\in G$. 
Given a subset $K\subseteq G\rtimes_\phi \Z$ and $r\in \Z$, we define $$K_r=\{x\in G\mid t^rx\in K\}=t^{-r}K\cap G.$$ 

In \cite{[BMMV06]}, the authors prove that [f.g. free]-by-cyclic groups have solvable conjugacy problem by reducing this question to the twisted conjugacy problem and Brinkmann's conjugacy problem on free groups. This was later generalized to more general extensions of groups in \cite{[BMV10]}, using orbit decidability, and very recently to ascending HNN-extensions of free groups in \cite{[Log22]} using variants of the $TCP$ and $BrCP$ for (nonsurjective) endomorphisms.  Similar ideas have also been explored in \cite{[CD22]} in the context of free-abelian times free groups, where it is proved that ascending HNN-extensions of free-abelian times free groups have solvable conjugacy problem.

In the same vein, the purpose of this paper is to relate the $GBrCP(G)$ and $GTCP(G)$ with $GCP(G\rtimes \Z)$.

The main result of the paper is the following:
\newtheorem*{generalizedmachine}{Theorem \ref{generalizedmachine}}
\begin{generalizedmachine}
Let $G$ be a finitely generated group, $\phi\in \Aut(G)$,  $K\subseteq G\rtimes_\phi \Z$ and $t^rg\in G\rtimes_\phi \Z$. Then:
\begin{enumerate}[i)]
\item if $r=0$, then $ GCP_{(K,t^rg)}(G\rtimes_\phi \Z)\equiv GBrCP_{(K_r,\phi,g)}(G)$;
\item if $r\neq 0$, then $GCP_{(K,t^rg)}(G\rtimes_\phi \Z) \equiv \bigplus\limits_{j=0}^{r-1} GTCP_{(K_r,\phi^r, g\phi^j)}(G)$. 
\end{enumerate}
\end{generalizedmachine}

Our main theorem is proved in a quite general form without imposing conditions on our target subsets and it provides us with an equivalence between an easier problem in $G\rtimes \Z$ and more complicated problems in $G$. However, as it will be made clear, even when a target set $K$ belongs to a well-behaved class of subsets, the subsets $K_r$ can be wild, which makes it difficult to apply one of the directions in some cases. 

We obtain corollaries from both directions of this equivalence: proving $GBrCP(G)$ and $GTCP(G)$ to solve $GCP(G\rtimes \Z)$ works better for recognizable and context-free subsets, while the converse works better for $f.g.\, cosets$, rational and algebraic subsets.

The main application concerns the class of virtually polycyclic groups. To do so, we prove the following theorem:

\newtheorem*{conjsep}{Theorem \ref{conjsep}}
\begin{conjsep}
Let $G$ be a conjugate separable finitely presented group such that $MP_{f.g.}(G)$  is decidable. Then $GCP_{[f.g. \, coset]}(G)$ is decidable.
\end{conjsep}

Using the connection between $GCP$ and $GBrCP$ and $GTCP$ that can be deduced from Theorem \ref{generalizedmachine}, we obtain the following corollary:
\newtheorem*{polyc br tcp}{Corollary \ref{polyc br tcp}}
\begin{polyc br tcp}
Let $G$ be a virtually polycyclic group. Then $GBrCP_{[f.g. coset]}(G)$ and $GTCP_{[f.g. coset]}(G)$ are decidable.
\end{polyc br tcp}

This is the most general result on the  $GBrCP$ to the knowledge of the author. Also, the solution of $GTCP$ also seems to be new. However, the simple version of the twisted conjugacy problem is solvable even for general endomorphisms of polycyclic groups by \cite{[Rom10]}.

Finally, we will discuss Brinkmann's equality problem: for the generalized version we will prove that decidability of rational intersection problem in $G\rtimes \Z$ yields decidability of $GBrP(G)$ and for the simple version we will prove that it is decidable for general endomorphisms of virtually free groups. 
\newtheorem*{brpvfree}{Theorem \ref{brpvfree}}
\begin{brpvfree}
Let $G$ be a f.g. virtually free group. Then $BrP_{End}$(G) is decidable.
\end{brpvfree}

The paper is organized as follows: in Section \ref{preliminaries}, we present some preliminaries on subsets of subgroups and on cyclic extensions. In Section \ref{genconjgbyz}, we prove the main result of the paper and deduce some natural corollaries for different classes of subsets: cosets of finitely generated groups, recognizable, context-free, rational and algebraic subsets. In Section \ref{vpoly} we apply the results obtained in the previous section to the class of virtually polycyclic groups, proving decidability of the coset versions of the generalized Brinkmann's  and twisted conjugacy problems. Finally, in Section \ref{brequality}, we discuss Brinkmann's equality problem, proving decidability of its simple version  for endomorphisms of virtually free groups and relating its rational generalized version with the rational intersection problem of cyclic extensions of the group.

\section{Preliminaries}
\label{preliminaries}
We will now present the basic definitions and results on rational, algebraic and context-free subsets of groups. For more detail, the reader is referred to \cite{[Ber79],[BS21],[Car22b]}. We will also discuss some notation and terminology and present a technical lemma concerning subsets of $G$-by-$\Z$ groups.

Let $G=\langle A\rangle$ be a finitely generated group, $A$ be a finite generating set, $\tilde A=A\cup A^{-1}$ and $\pi:\tilde A^*\to G$ be the canonical (surjective) homomorphism. 

A subset $K\subseteq G$ is said to be \emph{rational} if there is some rational language $L\subseteq \tilde A^*$ such that $L\pi=K$ and \emph{recognizable} if $K\pi^{-1}$ is rational.

We will denote by $Rat(G)$ and $Rec(G)$ the class of rational and recognizable subsets of $G$, respectively. Rational subsets generalize the notion of finitely generated subgroups.

\begin{theorem}[\cite{[Ber79]}, Theorem III.2.7]
\label{AnisimovSeifert}
Let $H$ be a subgroup of a group $G$. Then $H\in \text{Rat } G$ if and only if $H$ is finitely generated.
\end{theorem}

Similarly, recognizable subsets generalize the notion of finitely index subgroups.

\begin{proposition}[\cite{[Ber79]}, Exercise III.1.3]
\label{rec fi}
Let $H$ be a subgroup of a group $G$. Then $H\in \text{Rec } G$ if and only if $H$ has finite index in $G$.
\end{proposition}

In fact, if $G$ is a group and $K$ is a subset of $G$ then $K$ is recognizable if and only if $K$ is a (finite) union of cosets of a subgroup of finite index.

A natural generalization concerns the class of context-free languages. 
A subset $K\subseteq G$ is said to be \emph{algebraic} if there is some context-free language $L\subseteq \tilde A^*$ such that $L\pi=K$ and \emph{context-free} if $K\pi^{-1}$ is context-free. 
We will denote by $Alg(G)$ and $CF(G)$ the class of algebraic and context-free subsets of $G$, respectively.

It is obvious from the definitions that $Rec(G)$, $Rat(G)$, $CF(G)$ and $Alg(G)$ are closed under union, since both rational and context-free languages are closed under union. The intersection case is distinct: from the fact that rational languages are closed under intersection, it follows that $Rec(G)$ must be closed under intersection too. However, as it happens in the case of finitely generated subgroups, $Rat(G)$ is not necessarily closed under intersection. 
Another important closure property is given by the following lemma from \cite{[Her91]}.

\begin{lemma}
\label{cfletra}
\cite[Lemma 4.1]{[Her91]} Let $G$ be a finitely generated group, $R\in Rat(G)$ and $C \in \{Rat, Rec, Alg, CF\}$.
If $K \in C(G)$, then $KR, RK \in C(G).$
\end{lemma}

Let $A=\{a_1,\ldots, a_n\}$, $G=\langle A\mid R\rangle$ be a group and $\phi\in \Aut(G)$. A $G$-by-$\Z$ group admits a presentation of the form 
\begin{align*}
G\rtimes_\phi \Z=\langle A,t \mid R, t^{-1}a_it= a_i\phi\rangle.
\end{align*}

Every element of $G\rtimes_\phi \Z$ can be written in a unique way as an element of the form $t^ag$, where $a\in \Z$ and $g\in G$. 
Given a subset $K\subseteq G\rtimes_\phi \Z$ and $r\in \Z$, we define $$K_r=\{x\in G\mid t^rx\in K\}=t^{-r}K\cap G.$$

\begin{lemma}
\label{lem: cosets}
Let $G$ be a finitely generated group, $\phi\in \Aut(G)$, $H\leq_{f.g.} G\rtimes_\phi \Z$, $t^sg\in G\rtimes_\phi \Z$ and $K=(t^sg)H$. Then, for all $r\in \Z$, $K_r$ is either empty or a coset of $H\cap G$. Moreover, if we are given generators for $H$, we can decide whether $K_r$ is empty or  not and, in case it is nonempty, compute a coset representative for $K_r$.
\end{lemma}
\noindent\textit{Proof.} It is well known that the intersection of two cosets of a group is either empty or a coset of the intersection. Now, suppose that we are given generators for $H$, say $\{t^{k_1}g_1,\ldots , t^{k_n}g_n\}$. We have that
  $$K_r=\emptyset\Leftrightarrow (t^{s-r}g)H\cap G=\emptyset \Leftrightarrow H\cap t^{r-s} G=\emptyset \Leftrightarrow r-s\not\in \langle k_1 ,\ldots, k_n\rangle \leq \Z,$$
which is clearly decidable. Moreover, in case $K_r$ is nonempty, we can compute $\lambda_i\in \Z$ such that $r-s=\sum_{i=1}^n \lambda_i k_i$. Then, we compute  $h\in G$ such that  
$$t^rh=(t^sg)\prod_{i=1}^n(t^{k_i}g_i)^{\lambda_i}\in K\cap t^r G,$$ and so $h\in t^{-r}K\cap G=K_r$. Since $K_r$ is a coset, then $K_r= (K\cap G)h$.
\qed\\

 Notice that, if $K\leq_{f.g} G\rtimes_\phi \Z$, then, in general, $K_r=t^{-r}K\cap G$ is not a subgroup of $G$. In fact, it is a subgroup if and only if $t^r\in K$: 
   if $t^{r}\not\in K$, then $1\not\in t^{-r}K$ and so $1\not\in t^{-r}K\cap G$. If $t^r\in K$, then $t^{-r}K\cap G=K\cap G$, which is a subgroup of $G$.
 In particular, $K_0=K\cap G$ is always a subgroup, but not necessarily finitely generated. Obviously, if $G\rtimes_\phi\Z$ is Howson, then $K_0$ is finitely generated.

We will also write $G\rtimes \Z$ to denote the whole class of $G$-by-$\Z$ groups and so $GCP(G\rtimes \Z)$ will represent the uniform generalized conjugacy problem, i.e., taking the automorphism that defines the semidirect product as an input, while $GCP(G\rtimes_\phi \Z)$ will simply denote the generalized conjugacy problem for the group  $G\rtimes_\phi \Z$.

\section{Generalized conjugacy problem on $G$-by-$\Z$ groups}
\label{genconjgbyz}
The purpose of this section is to prove a generalized version of the result in \cite{[BMMV06]}, establishing a connection between $GCP(G\rtimes Z)$ and $GBrCP(G)$ and $GTCP(G)$ and discuss some possible applications for different classes of subsets.

\subsection{The main result}
In \cite{[BMMV06]}, the authors prove that [f.g. free]-by-cyclic groups have solvable conjugacy problem by reducing this to the twisted conjugacy problem and Brinkmann's conjugacy problem on free groups.  We now prove a generalized version of their result.

\begin{theorem}
\label{generalizedmachine}
Let $G$ be a group, $\phi\in \Aut(G)$,  $K\subseteq G\rtimes_\phi \Z$ and $t^rg\in G\rtimes_\phi \Z$. Then:
\begin{enumerate}[i)]
\item if $r=0$, then $ GCP_{(K,t^rg)}(G\rtimes_\phi \Z)\equiv GBrCP_{(K_r,\phi,g)}(G)$;
\item if $r\neq 0$, then $GCP_{(K,t^rg)}(G\rtimes_\phi \Z) \equiv \bigplus\limits_{j=0}^{r-1} GTCP_{(K_r,\phi^r, g\phi^j)}(G)$. 
\end{enumerate}
\end{theorem}
\noindent\textit{Proof.}
We start by proving i), so suppose that $r=0$. For $t^sv\in G\rtimes_\phi \Z$, we have that
 \begin{align*} 
  & v^{-1}t^{-s}t^0gt^sv\in K\\
  \Leftrightarrow  \quad &v^{-1}(g\phi^s)v\in K\\
    \Leftrightarrow  \quad &v^{-1}(g\phi^s)v\in K\cap G,
 \end{align*} 
and so $ GCP_{(K,t^0g)}(G\rtimes_\phi \Z)\equiv GBrCP_{(K_0,\phi,g)}(G)$.
 
To prove ii), let $r\neq 0$. For $t^sv\in G\rtimes_\phi \Z$, we have that
\begin{align*}
& v^{-1}t^{-s}t^rgt^sv\in K\\
\Leftrightarrow \quad & t^r(v^{-1}\phi^r)(g\phi^s)v\in K\\
\Leftrightarrow \quad & (v^{-1}\phi^r)(g\phi^s)v\in K_r\\
\end{align*}
Since 
 $$(v^{-1}\phi^r)(g\phi^s)v=((v^{-1}(g\phi^{s-r}))\phi^r)(g\phi^{s-r})(g^{-1}\phi^{s-r})v,$$ 
 then $g\phi^s$ has a $\phi^r$-twisted conjugate belonging to $K_r$ if and only if $g\phi^{s-r}$ has a $\phi^r$-twisted conjugate belonging to $K_r$. Hence, it suffices to check the existence of $\phi^r$-twisted conjugates for $0\leq s\leq r-1$, i.e. if there are $s\in \Z$ and $v\in G$ such that $(v^{-1}\phi^r)(g\phi^s)v\in K_r$ if and only if there are $0\leq s'\leq
r-1$ and $v'\in G$ such that  $({v'}^{-1}\phi^r)(g\phi^{s'})v'\in K_r$
\qed\\

This way, the study of the generalized conjugacy problem with respect to $K$ is related to the study of the generalized Brinkmann and twisted conjugacy problems with respect to $K_r$. Since, for $K\leq_{f.g} G\rtimes \Z$, $K_r$ is not necessarily a subgroup, the coset setting will be more adequate to us than the finitely generated subgroup setting, in view of Lemma \ref{lem: cosets}.

\subsection{The case of cosets}
 The following corollary is an immediate application of  Theorem \ref{generalizedmachine}  together with Lemma \ref{lem: cosets}.
\begin{corollary}
\label{corolario prob facil grupo dificil}
Let $G$ be a finitely generated group such that $G\rtimes \Z$ is Howson and $\phi\in \Aut(G)$. If for all $H\leq_{f.g} G\rtimes \Z$,  we can compute a finite set of generators for $H\cap G$, then:
\begin{enumerate}[i)]
\item if $GBrCP_{(f.g.,\phi)}(G)$ and $GTCP_{[f.g. coset]}(G)$ are decidable, then $GCP_{f.g}(G\rtimes_\phi \Z)$ is decidable;
\item if $GBrCP_{f.g.}(G)$ and $GTCP_{[f.g. coset]}(G)$ are decidable, then $GCP_{f.g}(G\rtimes \Z)$ is decidable;
\end{enumerate}
\end{corollary}

The following corollaries show us how we can do the converse. 
\begin{corollary}
\label{corolario maquina}
Let $G$ be a  group, $\phi\in \Aut(G)$,  $K\subseteq G$ and $g\in G$. Then:
\begin{enumerate}[i)]
\item $ GBrCP_{(K,\phi,g)}(G)  \equiv GCP_{( K,g)}(G\rtimes_\phi \Z)$;
\item $GTCP_{(K,\phi,g)}(G)  \equiv GCP_{(tK,tg)}(G\rtimes_\phi \Z)$.
\end{enumerate}s
\end{corollary}
\noindent\textit{Proof.} i) is an immediate application of Theorem \ref{generalizedmachine} i). To see ii), notice that if $K\subseteq G$ and $K'=tK\subseteq G\rtimes_{\phi}\Z$, then $$K_1'=t^{-1}K'\cap G=t^{-1}tK\cap G=K\cap G=K.$$ Now the result follows immediately from Theorem \ref{generalizedmachine} ii). 
\qed\\

\begin{corollary}
\label{corolario casos}
Let $G$ be a  finitely generated group  and $\phi\in \Aut(G)$.
 Then the following hold:
\begin{enumerate}[(i)]
\item if $GCP_{f.g}(G\rtimes_\phi \Z)$ is decidable, then $GBrCP_{(f.g., \phi)}(G)$ is decidable;
\item if $GCP_{f.g}(G\rtimes \Z)$ is decidable, then $GBrCP_{f.g.}(G)$ is decidable;
\item if $GCP_{[f.g.coset] }(G\rtimes_\phi \Z)$ is decidable, then $GBrCP_{([f.g.coset], \phi)}(G)$ and $GTCP_{([f.g.coset], \phi)}(G)$ are decidable;
\item if $GCP_{[f.g.coset]}(G\rtimes \Z)$ is decidable, then $GBrCP_{[f.g.coset]}(G)$ and $GTCP_{[f.g.coset]}(G)$ are decidable.
\end{enumerate}
\end{corollary}

Similar results hold for the \emph{simple} versions of these problems. For example, in 1969, it was proved in  \cite{[Rem69]} that polycyclic groups have solvable conjugacy problem. In particular $\Z^n\rtimes \Z$ has solvable conjugacy problem, and so $BrCP(\text{GL}_m(\Z))$ is decidable, which only became known in 1986 after Kannan and Lipton (\cite{[KL86]}) proved directly a (more general) version of Brinkmann's conjugacy problem for arbitrary matrices with rational entries.

\subsection{Other natural cases}

We saw how these problems relate when the class of subsets we consider is the class of cosets of finitely generated groups. In this case, one of the directions works \emph{better} than the other, in the sense that Corollary \ref{corolario prob facil grupo dificil} needs the (strong) additional hypothesis that $H\cap G$ is finitely generated and computable for all $H\leq_{f.g.} G$.

We will now see what we obtain when considering other natural classes of subsets of groups. In the rational and algebraic cases something similar will happen while in the recognizable and context-free cases we get the opposite as the analogous of Corollary \ref{corolario casos} will be harder to use.

\subsubsection{Recognizable subsets}

We start with a technical lemma.
\begin{lemma}
\label{Kr recognizable}
Let $G$ be a finitely generated group, $\phi\in \Aut(G)$ and $K\in Rec(G\rtimes_\phi \Z)$ . Then $K_r\in Rec(G)$, for all $r\in \Z$,
\end{lemma}
\noindent\textit{Proof.} Let $r\in \Z$. If $K\in Rec(G\rtimes_\phi \Z)$, then $t^{-r}K\in Rec(G\rtimes_\phi \Z)$ and so $t^{-r}K\cap G\in Rec(G)$ (\cite[III, Exercise 1.1]{[Ber79]}).  
\qed\\

The previous lemma, combined with Theorem \ref{generalizedmachine}, allows us to deduce that solving $GBrCP_{Rec}(G)$ and $GTCP_{Rec}(G)$ is enough to solve $GCP_{Rec}(G\rtimes \Z)$.
\begin{corollary}
\label{recimplications}
Let $G$ be a finitely generated group and $\phi\in \Aut(G)$. Then:
\begin{enumerate}[i)]
\item if $GBrCP_{(Rec,\phi)}(G)$ and $GTCP_{Rec}(G)$ are decidable, then $GCP_{Rec}(G\rtimes_\phi \Z)$ is decidable;
\item if $GBrCP_{Rec}(G)$ and $GTCP_{Rec}(G)$ are decidable, then $GCP_{Rec}(G\rtimes \Z)$ is decidable.
\end{enumerate}
\end{corollary}

The converse implication is not easy to use, since a recognizable subset of $G$ is not, in general, a recognizable subset of $G\rtimes \Z$ (for example $G\not\in Rec(G\rtimes \Z)$, since $[G\rtimes \Z:G]=\infty$).

We will now see another way of solving $GCP_{Rec}(G)$.
Recall that a subset $K$ of a group $G$ is recognizable if and only if it is a (finite) union of cosets of some finite index subgroup of $G$.
The finite index (resp. finitely generated) membership problem in a group $G$, $MP_{f.i.}(G)$ (resp. , $MP_{f.g.}(G)$), consists on deciding, given as input an element $g\in G$ and a finite index (resp. finitely generated) subgroup $H\leq G$, whether $g\in H$ or not. 

\begin{proposition}
If $MP_{f.i.}(G)$ is decidable, then $GCP_{Rec}(G)$ is decidable.
\end{proposition}
\noindent\textit{Proof.} Suppose that we can decide  $MP_{f.i.}(G)$. Let $K\in Rec(G)$ and $g\in G$. We want to decide if there is some $x\in G$ such that $x^{-1}gx\in K$. We know that $K$  is a (finite) disjoint union of cosets of some finite index subgroup $H\leq_{f.i.} G$, i.e.,
$$K=\bigcup_{i=1}^m Hb_i,$$
for some $m\in \N$,  $b_i\in K$. Also, $H$ has finitely many conjugates given by $b_i^{-1}Hb_i$ (and each of those is a finite index subgroup of $G$).

So, we only have to decide, for each $i\in [m]$ whether there is some $x\in G$ such that $x^{-1}gx\in Hb_i.$ But 
$$x^{-1}gx\in Hb_i\Leftrightarrow g\in xHb_ix^{-1} \Leftrightarrow g\in xHx^{-1}xb_ix^{-1}.$$

Since $H$ has finitely many conjugates, then there are finitely many possible values for $xHx^{-1}$ and since each of them has finite index in $G$, it also suffices to do finitely many checks on possible values $xb_ix^{-1}$, using $MP_{f.i.}(G)$. 
\qed\\

In \cite{[Har11]}, it is proved that every finitely $L$-presented group has solvable $MP_{f.i.}(G)$ and in \cite{[Rau22]}, it is proved that for recursively presented groups, $MP_{f.i.}(G)$ is equivalent to having computable finite quotients (CFQ). Hence, if $G$ is finitely $L$-presented or $G$ is a recursively presented group with CFQ, then $GCP_{Rec}(G)$ is decidable.

\subsubsection{Context-free subsets}
 
 We now prove the context-free analogous to Lemma \ref{Kr recognizable}.

\begin{lemma}
\label{Kr cfree}
Let $G$ be a finitely generated group, $\phi\in \Aut(G)$ and $K\in CF(G\rtimes_\phi \Z)$ . Then $K_r\in CF(G)$, for all $r\in \Z$,
\end{lemma}
\noindent\textit{Proof.}  If $K\in CF(G)$, by Lemma \ref{cfletra}, we have that $t^{-r} K\in CF(G\rtimes_\phi \Z)$. By \cite[Corollary 4.4 (b)]{[Her91]}, $K_r=t^{-r} K\cap G\in CF(G).$   
\qed\\

So we have the following corollary.
\begin{corollary}
\label{cfimplications}
Let $G$ be a finitely generated group  and $\phi\in \Aut(G)$. 
\begin{enumerate}[i)]
\item if $GBrCP_{(CF,\phi)}(G)$ and $GTCP_{CF}(G)$ are decidable, then $GCP_{CF}(G\rtimes_\phi \Z)$ is decidable;
\item if $GBrCP_{CF}(G)$ and $GTCP_{CF}(G)$ are decidable, then $GCP_{CF}(G\rtimes \Z)$ is decidable.
\end{enumerate}
\end{corollary}

Again, the converse implication is obstructed by the fact that a context-free subset of $G$ is not necessarily a context-free subset of $G\rtimes_\phi \Z$ (see \cite{[Car22b]}). 

\subsubsection{Rational and algebraic subsets}

It is well known that $F_2\rtimes_\phi \Z$ is not Howson, and so the intersection of rational subsets might not be rational. The next example shows something stronger:  $F_2\rtimes_\phi \Z$  has a rational subset $K$ such that $K_r$ is not a rational subset of $F_2$ for any $r\in \Z$.
\begin{example}
Take $F_2=\langle a,b\lvert\rangle$, $\phi=\lambda_b$. Then $K=(t\cup t^{-1})^*a(t\cup t^{-1})^*$ is a rational language of $\widetilde{\{a,b,t\}}^*$ and the subset of $F_2\rtimes_\phi \Z$ defined by $K$ is
$$\{t^{m+n}(a\phi^{n})\mid m,n\in \Z\}=\{t^{m+n}b^{-n}ab^n\mid m,n\in \Z\}.$$ 
We have that, for every $r\in \Z$,  $K_r=\{b^{-n}ab^{n}\in F_2\mid n\in \Z\}$ which is not a rational subset of $F_2$, by Benois' Theorem.
\end{example}
Hence, similarly to what happens in the coset case, in the rational case, we cannot expect $K_r$ to be a rational subset of $G$ when $K$ is a rational subset of $G\rtimes \Z$.

We remark that  if $K \in Rat(G\rtimes_\phi \Z)$, then, $K_r\in Alg(G\rtimes_\phi \Z)$, for all $r\in \Z$.  Indeed, let $A$ be a set of generators for $G$, $\pi:\widetilde{A\cup\{t\}}\to G\rtimes_\phi\Z$ be a surjective homomorphism and 
 $K\in Rat(G\rtimes_\phi \Z)$. Then $t^{-r}K\in Rat(G\rtimes_\phi \Z)$. Hence, there is a rational language $L\subseteq \widetilde{A\cup\{t\}}^*$ such that $L\pi=t^{-r}K$. The language $$L'=\{w\in L\mid |w|_t=|w|_{t^{-1}}\}=L\cap \{w\in \widetilde{ A\cup\{t\}}^*\mid |w|_t=|w|_{t^{-1}}\}$$ is context-free since it is the intersection of a rational language with a context-free language and $L'\pi=t^{-r}K\cap G$. Hence, $K_r\in Alg(G\rtimes_\phi \Z)$. 
However, we might have that $K_r\not\in Alg(G)$. In \cite[Example 4.7]{[Car22b]}, it is proved that if $G=\Z^2$, $\varphi=\begin{bmatrix}
    2 &1 \\
    1& 1\end{bmatrix}$ and $K$ is the orbit of $(1,0)$ through $\varphi$, then $K$ is rational in $\Z^2\rtimes_\varphi \Z$, but $K_0$ is not algebraic in $\Z^2$.  
So, even for rational $K$, $K_r$ might have a wild behaviour. For this reason, it is hard to obtain an analogous of Corollary \ref{corolario prob facil grupo dificil} in this case.

However, in the rational and algebraic case, unlike the recognizable and the context-free cases, the converse implication works without additional conditions. Indeed,  $K\in Rat(G)$  (resp.   $K\in Alg(G)$) implies that $K\in Rat(G\rtimes_\phi \Z)$ (resp. $K\in Alg(G\rtimes_\phi \Z)$) , for every $\phi\in \Aut(G)$ and so $xK\in  Rat(G\rtimes_\phi \Z)$ (resp.  $xK\in  Alg(G\rtimes_\phi \Z)$), for all $\phi\in \Aut(G)$ and $x\in G\rtimes_\varphi \Z$.

So we have the following corollary by directly applying Corollary \ref{corolario maquina}.
\begin{corollary}
\label{rat alg implications}
Let $G$ be a  finitely generated group  and $\phi\in \Aut(G)$. 
\begin{enumerate}[i)]
\item if $GCP_{Rat}(G\rtimes_\phi \Z)$ is decidable, then $GBrCP_{(Rat, \phi)}(G)$ and $GTCP_{(Rat, \phi)}(G)$ are   decidable;
\item if $GCP_{Rat}(G\rtimes \Z)$ is decidable, then $GBrCP_{Rat}(G)$  and $GTCP_{Rat}(G)$ are decidable;
\item if $GCP_{Alg}(G\rtimes_\phi \Z)$ is decidable, then $GBrCP_{(Alg, \phi)}(G)$ and $GTCP_{(Alg, \phi)}(G)$ are   decidable;
\item if $GCP_{Alg}(G\rtimes \Z)$ is decidable, then $GBrCP_{Alg}(G)$  and $GTCP_{Alg}(G)$ are decidable.
\end{enumerate}
\end{corollary}

An interesting consequence of Corollary \ref{rat alg implications} is that, in the case of  the context-free and recognizable versions of the conjugacy problem, decidability in $G$-by-$\Z$ groups yields the decidability in $G$-by-finite groups. To prove so, we start with a proposition that follows the strategy used in \cite{[LS11]} to prove the rational generalized conjugacy problem in virtually free groups. 

We denote by $\Via(G)$ the set of virtually inner automorphisms of a group $G$. 
\begin{proposition}
Let $\mc C\in\{Rec, Rat, CF, Alg\}$, $H$ be a finitely generated group and $G$ be an $H$-by-finite group. If $GTCP_{(\Via, \mc C)}(H)$ is decidable, then  $GCP_{\mc C}(G)$ is decidable.
\end{proposition}
\noindent\textit{Proof.} We have that $H$ is a finite index normal subgroup of $G$ and so $G$ admits a decomposition as a disjoint union
$$G= Hb_1\cup \ldots \cup Hb_m,$$ for some $b_i\in G$.
For each $i\in [m]$, we  define $\varphi_i:H\to H$ by $h\mapsto b_ihb_i^{-1}$. Then, $\varphi_i$ is clearly an automorphism of $H$ and since $\faktor{G}{H}$ is finite, then, for some $k\in \N$ we have that $b_i^k\in H$ and so $\varphi_i^k\in \Inn(H)$. Thus, $\varphi_i\in \Via(H)$.

Let $K\in \mc C(G)$, $h\in H$, $i\in [m]$ and assume that  $GTCP_{(Via, \mc C)}(H)$ is decidable. 
Then $hb_i$ has a conjugate in $K$ if and only if there are $h'\in H$ and $j\in [m]$ such that 
$(b_j^{-1}h'^{-1})(gb_i)(hb_j)\in K.$
Now,
\begin{align*}
&(b_j^{-1}h'^{-1})(gb_i)(hb_j)\in K\\
\Leftrightarrow \; &b_j^{-1}h'^{-1}g(h\varphi_i)b_ib_j\in K\\
\Leftrightarrow \; &(h'^{-1}g(h\varphi_i))\varphi_j^{-1}\in Kb_j^{-1}b_i^{-1}b_j \\
\Leftrightarrow \; &(h'^{-1}g(h\varphi_i))\varphi_j^{-1}\in (H\cap Kb_j^{-1}b_i^{-1}b_j) \\
\Leftrightarrow \;&h'^{-1}g(h\varphi_i)\in (H\cap Kb_j^{-1}b_i^{-1}b_j)\varphi_j. 
\end{align*}

Since $K\in \mc C(G)$, then $Kb_j^{-1}b_i^{-1}b_j\in \mc C(G)$ and $H\cap Kb_j^{-1}b_i^{-1}b_j\in \mc C(H)$ (see \cite{[Car22b]} for the cases of $CF$ and $Alg$ and \cite{[Ber79]} for the cases of $Rat$ and $Rec$) and so  $(H\cap Kb_j^{-1}b_i^{-1}b_j)\varphi_j\in \mc C(H).$

Therefore, our problem is now reduced to $m$ instances of $GTCP_{(\Via, \mc C)}(H)$, one for each $j\in [m]$, which we can solve. 
\qed\\

Combining the previous proposition with Corollary \ref{rat alg implications}, we obtain the following.
\begin{corollary}
Let $\mc C\in\{Rat, Alg\}$,  $G$ be a finitely generated group and $F$ be a finite group. If $GCP_{\mc C}(G\rtimes \Z)$ is decidable, then  $GCP_{\mc C}(G\rtimes F)$ is decidable.
\end{corollary}

 \section{Virtually polycyclic groups}
\label{vpoly}
A group $G$ is said to be \emph{conjugate separable} if, for all $g,h\in G$,

$$g \text{ and } h \text{ are conjugate }\Leftrightarrow g  \text{ and } h \text{ are conjugate in every finite quotient of $G$}.$$

If a group $G$ is finitely presented and conjugate separable, then it has a decidable conjugacy problem (see \cite{[Mos66]}). We now show that, adding the solution to the finitely generated membership problem, we can also solve the $[f.g. \, coset]$-generalized conjugacy problem. 

\begin{theorem}
\label{conjsep}
Let $G$ be a conjugate separable finitely presented group such that $MP_{f.g.}(G)$  is decidable. Then $GCP_{[f.g. \, coset]}(G)$ is decidable.
\end{theorem}
\noindent\textit{Proof.} Let $a,b\in G$ and $H=\langle h_1,\ldots, h_k\rangle \leq_{f.g.} G$ be our input. We want to decide if $a$ has a conjugate in $bH$. To do so, we will run two partial algorithms: one that stops,  answering \texttt{YES} if $a$ has a conjugate in $bH$ (but does not stop otherwise) and one that stops answering  \texttt{NO} if $a$ does not have a conjugate in $bH$ (but does not stop otherwise).

The first one is simple: we enumerate all conjugates of $a$ and check if they belong to $bH$ using $MP_{f.g.}(G)$. If we find a positive answer, we stop and answer  \texttt{YES}; while we don't we keep running.

Now we describe the second one. We start by enumerating all homomorphisms from $G$ onto finite groups. For each surjective homomorphism $\varphi:G\to F$, where $F=\{x_1,\ldots, x_m\}$ is a  finite group, we compute  $(bH)\varphi =b\varphi\langle h_1\varphi,\ldots, h_k\varphi \rangle $ and the set $X$ of all conjugates of $a\varphi$, which can be done since $F$ is finite. 
Compute $y_i\in G$ such that $y_i\varphi=x_i.$ We have that 

\begin{align*}
G=y_1 \Ker(\varphi)\cup \cdots \cup y_m\Ker(\varphi) \end{align*}
and 
$$bH=(y_1 \Ker(\varphi)\cap bH)\cup \cdots \cup (y_m\Ker(\varphi)\cap bH).$$
Each of the intersections $y_i\Ker(\varphi)\cap bH$ is either empty or a coset of $\Ker(\varphi)\cap H$. In fact, $y_i\Ker(\varphi)\cap bH\neq \emptyset$  if and only if $x_i\in (bH)\varphi$, which we can check. Moreover, for each nonempty intersection  $y_i\Ker(\varphi)\cap bH$, we can compute a representative $h_i\in bH$ such that  $y_i\Ker(\varphi)\cap bH=h_i(\Ker(\varphi)\cap H)$, by enumerating elements of $bH$ until we find one whose image is $x_i$.
Thus, we can write $bH$ as a disjoint union of the form 
$$bH=h_1(\Ker(\varphi)\cap H)\cup \cdots \cup h_n(\Ker(\varphi)\cap H),$$
where $n\leq m$ and $h_i\varphi=x_i$. Also, we can compute generators for $\Ker(\varphi)\cap H$ by Schreier's Lemma.

Now we check for each $h_i$ whether or not $h_i\varphi\in X$ and write $$bH\cap X\varphi^{-1}=\bigcup_{h_i\varphi\in X} h_i(\Ker(\varphi)\cap H).$$

This is a union of (computable) recognizable subsets of $H$, since $\Ker(\varphi)\cap H\leq_{f.i.} H$ and so $bH\cap X\varphi^{-1}$ is a (computable) recognizable subset of $H$. This gives us a set of candidate conjugators of $a$ in $H$ in the sense that it contains all elements that are mapped by $\varphi$ into a conjugator of $a\varphi$ (in $F$). 
So,  we keep computing the sets  $H\cap X\varphi^{-1}$, running through all homomorphisms $\varphi$ and intersecting them (which can be effectively done, since $Hb\cap X\varphi^{-1}\in Rec(H)$).  By conjugate separability, $a$ does not have a conjugate in $H$ if and only if at some point this intersection becomes $\emptyset$.
\qed\\

A group is \emph{polycyclic} if it admits a subnormal series 
$$G=G_0\triangleright G_1\triangleright\ldots G_n=\{1\}
$$
such that $\faktor{G_{i-1}}{G_i}$ is cyclic for $i\in [n]$.

In \cite{[Rem69]} and \cite{[For76]}, Remeslennikov and Formanek proved that virtually polycyclic groups are conjugate separable. Mal'cev in \cite{[Mal83]} proved that polycyclic groups are subgroup separable, and so, they have decidable membership problem, which implies that virtually polycyclic groups do too. So, we have the following corollary.
\begin{corollary}
\label{cor: gcp poly}
Let $G$ be a virtually polycyclic group. Then $GCP_{[f.g. coset]}(G)$ is decidable.
\end{corollary} 

If $G$ is polycyclic with subnormal series $$G=G_0\triangleright G_1\triangleright\ldots G_n=\{1\},
$$ such that $\faktor{G_{i-1}}{G_i}$ is cyclic for $i\in [n]$, 
and $\phi\in \Aut(G)$, then 
$$G\rtimes_\phi \Z \triangleright G_0 \triangleright G_1\triangleright\ldots G_n=\{1\}
$$
and $\faktor{G\rtimes_\phi \Z}{G_0}\simeq \Z$. Hence, $G\rtimes_\phi\Z$ is polycyclic.

We will now prove the same result for virtually polycyclic groups. 
A subgroup $H$ of a group $G$ is \emph{fully invariant} if $\phi(H)\subseteq H$ for every endomorphism $\phi$ of $G$.

\begin{proposition}
Let $G$ be a virtually polycyclic group and $\phi\in \Aut(G)$. Then $G\rtimes_\phi \Z$ is virtually polycyclic.
\end{proposition}
\noindent\textit{Proof.}
By \cite[Lemma 4.1]{[Car22]}, $G$ has a fully invariant finite index  normal polycyclic subgroup $H$: indeed if $P$ is a finite index polycyclic subgroup, then $P$ contains a finite index normal polycyclic subgroup $N$ and letting $H$ be the intersection of all normal polycyclic subgroups of index at most $[G:N]$ we obtain a fully invariant finite index normal polycyclic subgroup of $G$.

Since $H$ is fully invariant, the restriction $\psi=\phi|_H$ is an endomorphism of $H$. Write $G$ as a disjoint union 
$$G=H \cup Hb_1 \cup \cdots \cup Hb_n,$$
for some $b_n\in G$. Then 
$$G=G\phi= H\phi \cup  H\phi (b_1\phi) \cup \cdots \cup H\phi (b_n\phi),$$
and so $[G:H\phi]\leq n =[G:H]$. Since $[G:H\phi]=[G:H][H:H\phi]$, then $[G:H\phi]\geq [G:H]$. Hence, $[G:H\phi]=[G:H]$ and $[H:H\phi]=1$. 
Therefore, $\psi$ is bijective (injectivity is inherited from injectivity of $\phi$).

Since $H$ is polycyclic, then $H\rtimes_\psi \Z$ is polycyclic. Clearly $H\rtimes_\psi \Z\leq G\rtimes_\phi \Z$ and 
$$G\rtimes_\phi \Z= (H\rtimes_\psi \Z) \cup (H\rtimes_\psi \Z)t^0b_1 \cup \cdots \cup (H\rtimes_\psi \Z)t^0b_n,$$
and so   $G\rtimes_\phi \Z$ is virtually polycyclic.
\qed\\

From Corollaries \ref{corolario casos} and 
 \ref{cor: gcp poly} we deduce the following corollary.

\begin{corollary}
\label{polyc br tcp}
Let $G$ be a virtually polycyclic group. Then $GBrCP_{[f.g. coset]}(G)$ and $GTCP_{[f.g. coset]}(G)$ are decidable.
\end{corollary} 

Notice that Kannan and Lipton solved a version of Brinkmann's problem for matrices with entries in $\Q$. The generalized version of this algorithm is not easy: in \cite{[COW13]}, it is proved that it is decidable whether the orbit of a vector in $\Q^m$ by some matrix $A\in \Q^{m\times m}$ intersects a vector space of dimension at most $3$; for greater dimensions, the problem remains open. The case where matrices belong to GL$_m(\Z)$ and the target set is a coset of a f.g. subgroup is a consequence of Corollary \ref{polyc br tcp} 
and, to the author's knowledge, was unknown. 
Moreover, since in abelian groups $GBrCP$ and $GBrP$ coincide, we have the following corollary. We remark that, as highlighted in \cite{[Ven14]}, not much was known about $GBrP_{f.g}$ even for free-abelian groups. 

\begin{corollary}
Let $G$ be finitely generated abelian group. Then $GBrP_{[f.g. coset]}(G)$ is decidable.
\end{corollary}

\section{Brinkmann's equality problem}
\label{brequality}
Even though Brinkmann's conjugacy problem drew more attention than Brinkmann's equality problem, mostly due to its connection with the conjugacy problem established in \cite{[BMMV06]}, we believe that it is a question of independent interest: it is a very natural question concerning the dynamics of an automorphism and remains quite mysterious, specially in its generalized version, as mentioned by Ventura in \cite{[Ven14]}. Even in its simple version, not much is known: it is proved to be decidable for free-abelian groups in \cite{[KL86]}, for free groups in \cite{[Bri10]}, for braid groups in \cite{[GV14]} and for free-abelian-by-free groups in \cite{[CD22]}. Logan extended Brinkmann's result to general endomorphisms of the free group in \cite{[Log22]}.

In this section, we will relate the rational version of $GBrP(G)$ with the intersection problem of rational subsets of $G\rtimes \Z$ and solve  $GBrP_{(Rat,\Via)}(F)$ and  $GBrCP_{(Rat,\Via)}(F)$ for a free group $F$. For the simple version, we will prove it for virtually free groups. 

\subsection{Generalized version}

Naturally, one can reduce Brinkmann's problems restricted to $\Via$ to some conjugacy problems on the group. Using a result from \cite{[LS11]}, we can prove the following proposition. 

\begin{proposition}
\label{gbrpviafree}
Let $F$ be a finitely generated free group. Then $GBrP_{(Rat,\Via)}(F)$ and  $GBrCP_{(Rat,\Via)}(F)$ are decidable.
\end{proposition}
\noindent\textit{Proof.}
Let $\phi\in \Via(F)$ be such that $\phi^r=\lambda_x$, $u\in F$ and $K\in \Rat(F)$ be our input. Let $k\in \N$ and write $k=pr+q$ with $0\leq q\leq r-1$. Then $u\phi^k=(x^{-p}ux^p)\phi^q$. Thus, there is some $k\in \N$ such that $u\phi^k\in K$ (resp. $u\phi^k$ has a conjugate in $K$) if and only if there are $p \in \N$, $q\in\{0,\ldots, r-1\}$ such that $(x^{-p}ux^p)\phi^q\in K$ (resp. $(x^{-p}ux^p)\phi^q$ has a conjugate in $K$).

 For each $0\leq q\leq r-1$, the set $L=\{x^s\in F \mid s\in \N\}$ is rational and so $L\phi^q$ is rational too. The problem of deciding whether there is some $p\in \N$ such that $x^{-p}ux^p\in K$ is an instance of the rational generalized conjugacy problem with rational constraints, which is decidable for $F$ (see  \cite[Theorem 4.3]{[LS11]}).

 Moreover, deciding whether there is some $p\in \N$ such that $(x^{-p}ux^p)\phi^q$ has a conjugate in  $K$ can also be done, since this is equivalent to deciding whether $u\phi^q$ has a conjugate in $K$, which is  an instance of the rational generalized conjugacy problem
 \qed\\

The \emph{rational intersection problem}, $IP_{Rat}(G\rtimes \Z)$ consists on, given as input two rational subsets, deciding whether they intersect or not.
We now observe that, if we can solve $IP_{Rat}(G\rtimes \Z)$, then we can solve $GBrP_{Rat}(G)$.

\begin{proposition}
If $IP_{Rat}(G\rtimes \Z)$ is decidable, then $GBrP_{Rat}(G)$ is decidable.
\end{proposition}
\noindent\textit{Proof.}  Let $G$ be a group generated by a finite set $A=\{a_1,\ldots a_n\}$ and $\pi:\tilde A^*\to G$ be the standard surjective homomorphism. Then $G\rtimes_\phi \Z$ admits a presentation of the form (\ref{presentation semidirect}) and there is a natural surjective homomorphism $\rho:\widetilde{A\cup \{t\}}^*\to G\rtimes_\phi\Z$ such that $\rho|_A=\pi$ (identifying $G$ with the subset $\{t^0g\mid g\in G\}$). 
Let $K\in Rat(G)$, $\phi\in \Aut(G)$ and $g\in G$.  We want to decide if there is some $n\in \N$ such that $g\phi^n\in K$. Since $K$ is a rational subset of $G$, then it is a rational subset of $G\rtimes_\phi \Z$, and so is $t^*K.$ Let $w\in \tilde A^*$ be a word such that $w\pi=g$ and put $L=\{wt^n\mid n\in \N\}\subseteq \widetilde{A\cup \{t\}}^*$. Then $L$ is a rational language and $L\pi=\{t^ng\phi^n\mid n\in \N\}$ is a rational subset of $G\rtimes_\phi Z$. Now, there exists some $n\in \N$ such that 
$g\phi^n\in K$ if and only if $L\pi\cap t^*K\neq \emptyset$, which we can decide by hypothesis.
\qed\\

\subsection{Simple version}
We will now prove Brinkmann's equality problem for endomorphisms of finitely generated virtually free groups.

When we take a finitely generated virtually free group as input, we assume that we are given a decomposition as a disjoint union
\begin{align}
\label{decomp}
G=Fb_1\cup Fb_2\cup \cdots \cup Fb_m,
\end{align}
where $F=F_A\trianglelefteq G$ is a finitely generated fully invariant free group and a presentation of the form $\langle A,b_1,\ldots, b_m \mid R\rangle$, where the relations in $R$ are of the  form $b_ia=u_{ia}b_{i}$ and $b_{i}b_{j}=v_{ij}b_{r_{ij}}$, with  $u_{ia}, v_{ij} \in F_A$ and $r_{ij}\in [m]$, $i,j=1\ldots,m$, $a\in A$.

Given a group $G$, an element $g\in G$ and an endomorphism $\varphi\in \End(G)$, we denote by $\Orb_\varphi(g)$, the $\varphi$-orbit of $g$:
$$\Orb_\varphi(g)=\{g\varphi^k\mid k\in \N\}.$$

\begin{theorem}
\label{brpvfree}
Let $G$ be a finitely generated virtually free group. Then $BrP_{End}$(G) is decidable.
\end{theorem}
\noindent\textit{Proof.} Let $\varphi\in \End(G)$ and $g,h\in G$ be our input. We can assume that $g,h$ are given as $ub_i$ and $vb_j$ where $u,v\in F$ and $i,j\in[m]$. So, we want to decide if there is some $k\in \N$ such that $(ub_i)\varphi^k=vb_j$. Define the mapping $\theta: G/F\to G/F$ by $[b_i]\mapsto [b_i\varphi]$. The mapping is well defined since $F$ is fully invariant: given $w,w'\in F$, then, for all $i\in[m]$, we have that $$[(wb_i)\varphi]=[w\varphi b_i\varphi]=[b_i\varphi]=[w'\varphi b_i\varphi]=[(w'b_i)\varphi].$$ Moreover, it is an endomorphism of $G/F$: for all $i,j\in[m]$, 
$$([b_i]\theta)([b_j]\theta)=[b_i\varphi][b_j\varphi]=[(b_ib_j\varphi)]=[b_ib_j]\theta=([b_i][b_j])\theta.$$
Since $G/F$ is finite, we can compute the entire $\theta$-orbit of $[b_i]$, which must be finite.
Hence, we can verify if $[b_j]\in \Orb_\theta ([b_i])$, and the set $\{r\in \N\mid [b_i]\theta^r=[b_j]\}$ is either empty or of the form $s+p\N$ with  computable $s,p\in \N$.

If, for some $k\in \N$, $(ub_i)\varphi^k=vb_j$, then $vb_j=(u\varphi^k)(b_i\varphi^k)\in [b_i\varphi^k]=[b_i]\theta^k$. So, our candidate values for $k$ are precisely $s+p\N$. Indeed, we denote by $s$ the first time that $[b_j]$ occurs in the $\theta$-orbit of $[b_i]$ and by $p$, the $\theta$-period of $[b_j]$ (which might be $0$, in which case there is only one candidate).  Compute $y\in F$ such that $(ub_i)\varphi^s=yb_j$ and  $z\in F$  such that $b_j\varphi^p=zb_j$. It is easy to check by induction that, for all $d\in \N$, we have $$b_j\varphi^{dp}=\left(\prod\limits_{i=0}^{d-1}z\varphi^{ip}\right)b_j.$$

If $y=v$ then we answer yes and output $s$. Otherwise, let $c$ be a new letter, not belonging to $A$ and let $\psi:F*\langle c|\rangle\to F*\langle c|\rangle$ be defined by mapping $x$ to $x\varphi^p$ for every $x\in F$ and $c$ to $zc$. Also by induction, we can check that $$c\psi^d=\left(\prod\limits_{i=0}^{d-1}z\varphi^{ip}\right)c,$$ for all $d\in \N$.

 We claim that there is some $k\in \N$ such that $(ub_i)\varphi^k=vb_j$ if and only if there is some $k\in \N$ such that $(yc)\psi^k=vc$, which can be decided using BrP($F_n$). We have that there is  $k\in \N$ such that $(ub_i)\varphi^k=vb_j$ if and only if there is some $d\in \N$ such that  $(ub_i)\varphi^{s+pd}=vb_j$, i.e., if there is some $d\in \N$ such that  $(yb_j)\varphi^{pd}=vb_j$.
 We then have the following series of implications, which concludes the proof:
 \begin{align*}
& vb_j=  (yb_j)\varphi^{pd}\\
 \Leftrightarrow\quad &vb_j=(y\varphi^{pd})(b_j\varphi^{pd})\\
  \Leftrightarrow \quad&vb_j=(y\varphi^{pd})\left(\prod\limits_{i=0}^{d-1}z\varphi^{ip}\right)b_j\\
   \Leftrightarrow \quad&v=(y\varphi^{pd})\left(\prod\limits_{i=0}^{d-1}z\varphi^{ip}\right)\\
    \Leftrightarrow\quad &v=(y\psi^d)\left(\prod\limits_{i=0}^{d-1}z\varphi^{ip}\right)\\
      \Leftrightarrow\quad &vc=(y\psi^d)(c\psi^d)\\
          \Leftrightarrow \quad&vc=(yc)\psi^d.
          \end{align*}
\qed\\

We do not know whether the Brinkmann's conjugacy problem is decidable for virtually free groups or not. If it is, then that would imply that $CP(G\rtimes \Z)$ would be decidable
whenever $G$ is a finitely generated virtually free group, as $TCP_{Aut}(G)$ is decidable in that case \cite[Theorem 4.8]{[BMV10]}. In fact, knowing that the computation of  fixed subgroups for endomorphisms of virtually free groups is possible by \cite[Theorem 4.2]{[Car22]}, we can easily extend this result to general endormorphisms of finitely generated virtually free groups. 
\begin{theorem}
Let $G$ be a finitely generated virtually free group and $\varphi\in \End(G)$ be given algorithmically. Then $TCP_{End}(G)$ is decidable.
\end{theorem}
\noindent\textit{Proof.}
Let $\varphi\in \End(G)$ and $u,v\in G$ be our input.
Let $x,y$ be new letters not belonging to the generators of $G$ and put $G_1=G*\langle x,y|\rangle$. Then $G_1$ is virtually free, since it is a free product of virtually free groups. Let $\psi:G_1\to G_1$ be defined by $g\psi=g\varphi$ for all $g\in G$; $x\psi=xv$; and $y\psi=u^{-1}y$. Then, there exists a twisted conjugator $w\in G$ such that $v=w^{-1}u(w\varphi)$ if and only if there exists a fixed point of $\psi$ of the form $xzy$ with $z\in G$. Indeed, if $v=w^{-1}u(w\varphi)$, then 
$$(xw^{-1}y)\psi=xv(w^{-1}\varphi) u^{-1}y=xw^{-1}y$$ and, conversely, if there is a fixed point $xzy$ of $\psi$ with $z\in G$, then 
$$xzy=(xzy)\psi=xv(z\varphi) u^{-1}y,$$
hence $z=v(z\varphi) u^{-1}$.

So, $u$ and $v$ are $\varphi$-twisted conjugates if and only if $x^{-1}\Fix(\psi)y^{-1}\cap G\neq \emptyset$, which can be decided, since $G$ and $\Fix(\psi)$ are computable subgroups of $G_1$.
\qed\\

\section*{Acknowledgements}
The author is grateful to Pedro Silva for fruitful discussions of these topics, which  improved the paper.
The author was  supported by the grant SFRH/BD/145313/2019 funded by Funda\c c\~ao para a Ci\^encia e a Tecnologia (FCT).%

\bibliographystyle{plain}
\bibliography{Bibliografia}

 \end{document}